\documentclass[10pt, a4paper]{article}
\usepackage{graphics}
\usepackage{latexsym}
\usepackage{amssymb}
\usepackage{amsmath}
\usepackage{color}
\usepackage{comment}
\usepackage[title]{appendix}

\newcommand{\beq}{\begin{equation}}
\newcommand{\eeq}{\end{equation}}
\newcommand{\beqas}{\begin{eqnarray*}}
\newcommand{\eeqas}{\end{eqnarray*}}
\newcommand{\ep}{\varepsilon}
\newcommand{\ue}{u^{\varepsilon}}

\newcommand{\wc}{\rightharpoonup}

\newcommand{\bV}{{\bf V}}

\newcommand{\be}{\begin{equation}}
\newcommand{\ee}{\end{equation}}
\newcommand{\bea}{\begin{eqnarray}}
\newcommand{\eea}{\end{eqnarray}}
\newcommand{\beas}{\begin{eqnarray*}}
\newcommand{\eeas}{\end{eqnarray*}}

\def\IN{\mathbb{N}}

\def\IP{\mathbb{P}}
\def\IR{\mathbb{R}}
\def\IZ{\mathbb{Z}}
\def\IE{\mathbb{E}}

\newcommand{\vertiii}[1]{{\left\vert\kern-0.25ex\left\vert\kern-0.25ex\left\vert #1 
    \right\vert\kern-0.25ex\right\vert\kern-0.25ex\right\vert}}
\DeclareMathOperator*{\esssup}{ess\,sup}

\newtheorem{theorem}{Theorem}[section]
\newtheorem{lemma}[theorem]{Lemma}

\newtheorem{proposition}[theorem]{Proposition}
\newtheorem{definition}[theorem]{Definition}
\newtheorem{remark}[theorem]{Remark}

\numberwithin{equation}{section}
\begin{document}
\title{Optimal finite elements for  ergodic stochastic two-scale elliptic equations}
\author{Viet Ha Hoang, Chen Hui Pang and Wee Chin Tan\\ \\
Division of Mathematical Sciences\\ School of Physical and Mathematical Sciences\\ Nanyang Technological University, 21 Nanyang Link, Singapore 637371}
\date{}
\maketitle
\begin{abstract}
We develop an essentially optimal finite element approach for solving ergodic stochastic two-scale elliptic equations whose two-scale coefficient may depend also on the slow variable. We solve the limiting stochastic two-scale homogenized equation obtained from the stochastic two-scale convergence in the mean (A. Bourgeat, A. Mikelic and S. Wright, J. reine angew. Math, Vol. 456, 1994), whose solution comprises of  the solution to the homogenized equation and the corrector,  by truncating the infinite domain of the fast variable and using the sparse tensor product finite elements. We show that the convergence rate in terms of the truncation level is equivalent to that for solving the cell problems in the same truncated domain. Solving this equation, we obtain the solution to the homogenized equation and the corrector at the same time, using only a number of degrees of freedom that is essentially equivalent to that required for solving one cell problem. Optimal complexity is obtained when the corrector possesses sufficient regularity with respect to both the fast and the slow variables. Although the regularity norm of the corrector depends on the size of the truncated domain, we show that the convergence rate of the approximation for the solution to the homogenized equation is independent of the size of the truncated domain. With the availability of an analytic corrector, we construct a numerical corrector for the solution of the original stochastic two-scale equation from the finite element solution to the truncated stochastic two-scale homogenized equation. Numerical examples of quasi-periodic two-scale equations, and a stochastic two-scale equation of the checker board type, whose coefficient is discontinuous, confirm the theoretical results. 
\end{abstract}

\section{Introduction}
We develop an essentially optimal finite element (FE) method to solve ergodic stochastic two-scale problems. For these problems, the homogenized coefficient can be found from solutions of abstract cell problems that are posed in the abstract probability space (see, e.g, \cite{BLP}, \cite{JKO}, \cite{Armstrongbook}, \cite{GloriaOtto}). Each corresponding real realization of the solution of these abstract equations satisfies an equation in the whole real space. To establish the homogenized coefficient numerically, we need to solve cell problems in a truncated real domain. The accuracy of the approximation of the homogenized coefficient, which is also the accuracy of the approximation of the solution to the homogenized equation, depends essentially on the size of the truncated domain (\cite{BourgeatPiatnitski}). When the two-scale coefficient depends also on the slow variable, homogenization is more complicated to justify (\cite{BourgeatMW}). Further, for each macroscopic point in the slow variable domain, we need to solve  a set of cell problems to approximate the homogenized coefficient at that point, leading to a large level of complexity.

	Motivated by our previous works of solving locally periodic multiscale problems, we develop in this paper a high dimensional FE approach that provides the solution of the homogenized equation that approximates the solution of the two-scale equation macroscopically, and the corrector that encodes the microscopic information, at the same time. For the case of  locally periodic problems, the multiscale homogenized equation (\cite{Nguetseng, Allaire}) contains all the necessary information on the solution of the multiscale problem: its solution provides the solution to the homogenized equation and the corrector. However, this equation is posed in a high dimensional tensorized domain. Schwab \cite{Schwabicm} and Hoang and Schwab \cite{HSelliptic} develop the sparse tensor product FE approach that solves this equation with an essentially optimal level of complexity for obtaining an approximation of its solution within a prescribed level of accuracy. The method has been successfully applied to other classes of equations: wave equation (\cite{XHwave}), elasticity (\cite{XHelasticity, XHelasticwave}), electro-magnetic equations (\cite{ChuHoangMaxwelltype, ChuHoangMaxwellwave}), nonlinear equations
 (\cite{Hoangmonotone, TanHoangJCAMHighdimensional2019, TanHoangIMAJNASparsetensor2020}). For ergodic stochastic two-scale equations, the corresponding stochastic version of the two-scale homogenized equation is obtained from the stochastic two-scale convergence in the mean which is developed by
 Bourgeat et al. \cite{BourgeatMW}. In the stochastic two-scale homogenized equation  \eqref{eq:2shomprob},  $u_0(x)$ and $U_1(x,\omega)$ are the solution of the homogenized equation and the abstract corrector which is obtained from $u_0$ and the solutions of the abstract cell problems respectively. We show that $u_0$ and the realization $U_1(x,T(y)\omega)$ of the corrector (here $T(y)$ is the ergodic dynamical system on the probability space) can be approximated by the solution of the truncated stochastic two-scale homogenized equation which depends on real variables $x$ and $y$ where $x$ belongs to the real physical domain, and $y$ belongs to a truncated domain. Solving this equation with the sparse tensor product FEs,  we obtain an approximation whose accuracy in terms of the FE mesh width and the truncation level of the domain of the fast variable $y$ is essentially equal to the accuracy of the approximation of the homogenized coefficient, obtained by solving the cell problems in the same truncated domain using the same FE mesh width. We thus obtain an approximation for the solution of the homogenized equation and the corrector at the same time, using an essentially optimal number of degrees of freedom which is equivalent (apart from a possible logarithmic multiplying factor) to that required for solving one cell problem in the truncated domain, for the same level of accuracy. The essential optimality of the sparse tensor product FEs is achieved when the corrector possesses sufficient regularity. However, in this case, the regularity norm of the corrector grows with the size of the truncated domain. Nevertheless, we show that the approximation rate for the solution $u_0$ of the homogenized equation  is independent of the size of the truncated domain. Given an analytic corrector, we construct a numerical corrector from the numerical solution of the truncated stochastic two-scale  homogenized equation. Numerical experiments confirm our theory. 

We note that general multiscale methods have been applied for ergodic stochastic multiscale problems. We mention exemplarily the work by Efendiev and Pankov \cite{EfendievPankov} using the Multiscale Finite Element method (\cite{HouWu, EfendievHou}) and Gallistl and Peterseim \cite{GallistlPeterseim} using the localization approach (\cite{MalqvistPeterseim}). These approaches need to solve multiscale local problems using microscopic meshes that must be at most of the order of the microscopic scale. In this paper, our mesh width is macroscopic; the convergence rate depends on the level of the truncation of the fast variable domain, which depends solely on the statistics of the ergodic random field, and not on the microscopic scale.

The paper is organized as follows.

In the next section, we formulate the ergodic two-scale equation and review some results on their homogenization. In section   \ref{sec:twoscaleconv}, we recall the concept of stochastic two-scale convergence in the mean of Bourgeat et al. \cite{BourgeatMW}. We then establish the approximation of this equation by the truncated stochastic two-scale homogenized equation, where we consider a truncated domain for the fast variable with size $2\rho$, and show that the approximation accuracy in terms of the truncated level of the fast variable domain is equal to that of the cell problems. In Section \ref{sec:FE}, we solve the truncated stochastic two-scale homogenized problem by tensor product FEs. We show that the sparse tensor product FEs obtain essentially equal convergence rate as the full tensor product FEs using a far  smaller number of degrees of freedom, which is essentially optimal, given that the solution possesses sufficient regularity.  In Section \ref{sec:corrector}, with the availability of an analytic corrector result, we construct a numerial corrector for the solution of the two-scale equation from the FE solution of the truncated stochastic two-scale homogenized equation, when the microscopic scale, the FE mesh width converge to 0, and the size of the truncated domain tends to infinity. Section \ref{sec:numerical} presents some numerical experiments. We solve two-scale quasi periodic problems in one and two dimensional domains where the exact homogenized coefficient can be approximated highly accurately, so a reference solution can be approximated with a high level of accuracy. The numerical results confirm the theory. When the mesh size is larger, the decaying rate of the total error is similar to that of the FE error as the effect of  truncation of the fast variable domain is less significant. However, when the FE mesh gets smaller, the decaying rate deteriorates due to the more prominent effect of the truncated fast variable domain. When we increase the size of the truncated domain, the total error gets smaller. We then solve a two dimensional checker board type problem whose coefficient is piecewise constant. We observe similar behaviour in the total error. The appendices contain the involved proofs of  some technical results. 

Throughout the paper, by $c$ we denotes a generic constant whose value may change between different appearances.  Repeated indices indicate summation.
\section{Problem setting}\label{sec:setting}
\subsection{Ergodic two-scale problem}
Let $D$ be a bounded Lipschitz domain in $\IR^d$. Let $(\Omega,\Sigma, \IP)$ be a probability space. We assume that there is an ergodic dynamical system $T:\IR^d\times \Omega\to \Omega$ such that for $y, y'\in \IR^d$ and $\omega\in\Omega$, $T(y+y')\omega=T(y)T(y')\omega$. Furthermore, any invariant subsets of $\Omega$ have either probability 0 or 1. Let $A: \IR^d\times \Omega \to \IR^{d\times d}_{sym}$. There are positive constants $c_1$ and $c_2$ such that for all $\xi,\zeta\in\IR^d$
\be
c_1|\xi|^2\le A(x,\omega)\xi\cdot\xi,\ \ \mbox{and}\ \ A(x,\omega)\xi\cdot\zeta\le c_2|\xi||\zeta|,
\label{eq:boundednesscoerciveness}
\ee
for a.a. $x\in D, \omega\in\Omega$.
 We assume further that $A\in C(\bar D, L^{\infty}(\Omega))^{d\times d}$. Let $\ep>0$ be a small quantity that represents the microscopic scale that the problem depends on. We define the ergodic random two-scale coefficient as 
\[
A^\ep(x,\omega)=A(x,T({x\over\ep})\omega).
\]
We denote by $V=H^1_0(D)$ and $H=L^2(D)$.
Let $f\in V'$.  We consider the problem
\be
-\nabla\cdot(A^\ep(x,\omega)\nabla\ue)=f 
\label{eq:2sprob}
\ee
with the Dirichlet boundary condition $\ue(x)=0$ for $x\in\partial D$. 
\subsection{Homogenized equation}
We review in this section the well known results on homogenization of problem \eqref{eq:2sprob}. We refer to the standard references such as \cite{BLP} and \cite{JKO} for details. We denote by $\partial=(\partial_1,\ldots,\partial_d)$ where $\partial_i$ is the generator of the semingroup $T(x+e_i)$ and $e_i$ is the $i$th unit vector in the standard basis of $\IR^d$. 
Following \cite{JKO} and \cite{BourgeatPiatnitski}, we define $L^2_{pot}(\Omega)$ the completion of the space of all functions in $L^2(\Omega)^d$ of the form $\partial\psi$ for $\psi$ belonging to the domain of $\partial$; and by $L^2_{sol}(\Omega)$ the completion in $L^2(\Omega)^d$ of all functions $\psi$ whose components belong to the domain of the $\partial$ operator such as ${\rm div}\psi:=\sum_{i=1}^d\partial_i\psi=0$.  For each basis vector $e^i$ in the standard basis of $\IR^d$, we consider the cell problem: Find $w^i(x,\omega)\in L^2_{pot}(\Omega)$ such that
\be
A(x,\omega)(w^i(x,\omega)+e^i)\in L^2_{sol}(\Omega).
\label{eq:abstractcell}
\ee
Then the homogenized coefficient can be expressed as 
\be
A^0_{ij}(x)=\int_\Omega A_{ik}(x,\omega)(w^j_k(x,\omega)+\delta_{kj})\IP(d\omega).
\label{eq:A0}
\ee
We have that when $\ep\to 0$, $\ue(\cdot,\omega)\wc u_0$ in $V$ a.s. where the deterministic function $u_0$ satisfies the homogenized problem
\be
-\nabla\cdot(A^0(x)\nabla u_0)=f
\label{eq:homprob}
\ee
with the Dirichlet boundary condition on $\partial D$. We note that Jikov et al. \cite{JKO} only consider the case where $A(x,\omega)=A(\omega)$ i.e. it does not depend on $x$, but the results hold identically for the case where $A$ depends on $x$ as considered in this paper. 
For $\rho>0$, let $Y^\rho$ be the cube $(-\rho,\rho)^d\subset\IR^d$. We consider the problem
\be
-\nabla_y\cdot(A(x,T(y)\omega))(\nabla_y N^{i\rho}(x,y,\omega)+e^i)=0,\ \ y\in Y^\rho
\label{eq:cellprob}
\ee
with the Dirichet boundary condition $N^{i\rho}(x,y,\omega)=0$ for $y\in\partial Y^\rho$. We define the coefficient $A^{0\rho}(x,\omega)\in \IR^{d\times d}_{sym}$ as
\be
A_{ij}^{0\rho}(x,\omega)={1\over |Y^\rho|}\int_{Y^\rho}A_{ik}(x,T(y)\omega)({\partial N^{j\rho}\over\partial y_k}(x,y,\omega)+\delta_{kj})d y.
\label{eq:A0rho}
\ee
We then have
\begin{proposition}\label{prop:A0conv}
Assume that $\esssup_{\omega\in\Omega}|A(x,\omega)-A(x',\omega)|\le c|x-x'|$ for all $x,x'\in D$, where $c$ is independent of $x$ and $x'$. Almost surely, $\lim_{\rho\to \infty}\|A^0(\cdot)-A^{0\rho}(\cdot,\omega)\|_{L^\infty(D)^{d\times d}}=0$.
\end{proposition}
Bourgeat and Piatnitski \cite{BourgeatPiatnitski} show this result in the case where $A$ does not depend on $x$, i.e. $A=A(\omega)$ (see also \cite{GloriaOtto, Armstrongbook}).
We show this proposition for the case where $A$ also depends on $x$ in Appendix \ref{app:B}. 
We note that with further assumption on the structure of the probability space $(\Omega,\Sigma,\IP)$, we also have the stronger convergence result (\cite{BourgeatPiatnitski}):
%
\be
\IE[\|A^0-A^{0\rho}\|_{L^\infty(D)^{d\times d}}^2]\le c\rho^{-\beta}
\label{eq:A0rhoconvergence}
\ee
for $\beta>0$.
%
\section{Stochastic two-scale convergence in the mean}
\label{sec:twoscaleconv}
In this section we recall the stochastic two-scale convergence in the mean and the stochastic two-scale homogenized equation obtained from it of Bourgeat et al. \cite{BourgeatMW}. We then study approximation of its solution by the truncated stochastic two-scale homogenized equation. 
\subsection{Stochastic two-scale convergence in the mean}
The homogenized problem can also be established by using the "stochastic two-scale convergence in the mean" as defined by Bourgeat et al. \cite{BourgeatMW}. This extends the concept of two-scale homogenization (\cite{Nguetseng},\cite{Allaire}) to the stochastic setting. Following \cite{BourgeatMW} we say that a function $\psi\in L^2(D\times\Omega)$ is admissible if $\psi(x,T(x)\omega)$ belongs to $L^2(D\times\Omega)$. We first recall the definition
\begin{definition}
A sequence $\{w^\ep\}_\ep$ in $L^2(D\times\Omega)$  stochastically two-scale converges in the mean to a function $w\in L^2(D\times\Omega)$ if for every admissible function $\psi\in L^2(D\times\Omega)$,
\[
\lim_{\ep\to 0}\int_\Omega\int_D w^\ep(x,\omega)\psi(x,T({x\over\ep})\omega)dx\IP(d\omega)= \int_\Omega\int_D w(x,\omega)\psi(x,\omega)dx\IP(d\omega).
\]
\end{definition}
This definition makes sense due to the following result (\cite{BourgeatMW} Theorem 3.4):
\begin{lemma}
From a bounded sequence in $L^2(D\times\Omega)$ we can extract a stochastic two-scale convergent in the mean subsequence.
\end{lemma}
For a bounded sequence in $L^2(\Omega,V)$ we have the following result (\cite{BourgeatMW} Theorem 3.7):
\begin{lemma}
Let $\{w^\ep\}_\ep$ be a bounded sequence in $L^2(\Omega,V)$. There is a subsequence (not renumbered), a function $w\in V$ and a function $w_1\in L^2(D, L^2_{pot}(\Omega))$ such that
$w^\ep$ stochastically two-scale converges in the mean to $w$ and $\nabla w^\ep$ stochastically two-scale converges in the mean to $\nabla w+w_1$.
\end{lemma}
Using these, we  have the following result on the stochastic two-scale convergence in the mean limit for the solution of problem \eqref{eq:2sprob}. 
\begin{proposition}
There are functions $u_0\in V$ and $U_1\in L^2(D, L^2_{pot}(\Omega))$ such that the sequence $\{\ue\}_\ep$ of solution to problem \eqref{eq:2sprob} stochastic two-scale converges in the mean to $u_0$, and $\{\nabla\ue\}_\ep$ stochastic two-scale converges in the mean to $\nabla u_0+U_1$. The functions $u_0$ and $U_1$ form the unique solution of the problem:
\be
\int_\Omega\int_DA(x,\omega)(\nabla u_0(x)+U_1(x,\omega))\cdot(\nabla\phi_0(x)+\Phi_1(x,\omega))dx\IP(d\omega)=\int_Df(x)\phi_0 (x)dx,
\label{eq:2shomprob}
\ee
$\forall\,\phi_0\in V,$ $\Phi_1\in L^2(D,L^2_{pot}(\Omega))$.
\end{proposition}
We refer to Bourgeat et al. \cite{BourgeatMW} Theorem 4.1.1 for a proof.
\subsection{The truncated equation}
 We note the following result.
\begin{lemma}
Let $F_i\in L^2(\Omega)$ ($i=1,\ldots,d$), and  $G\in L^2(\Omega)$ be such that
\[
\int_\Omega(F_i(\omega)\partial _i\psi(\omega)+G(\omega)\psi(\omega))\IP(d\omega)=0
\]
for all $\psi$ belonging to the domain of the generator $\partial_i$ of the semi group $T(y+e_i)$ for all $i=1,\ldots,d$. Then a.s. for all $\phi\in C^\infty_0(\IR^d)$
\[
\int_{\IR^d}\left(F_i(T(y)\omega){\partial\phi\over\partial y_i}(y)+G(T(y)\omega)\phi(y)\right)dy=0.
\]
\end{lemma}
The proof of this result can be found, e.g. in Beliaev and Kozlov \cite{BeliaevKozlov} page 12 (though Beliaev and Kozlov consider the more complicated case of a randomly perforated domain). From \eqref{eq:2shomprob}, we have
\[
\int_\Omega\int_DA(x,\omega)(\nabla u_0(x)+U_1(x,\omega))\cdot(\nabla\phi_0(x)+\psi(x)\partial\phi_1(\omega))dx\IP(d\omega)=\int_Df(x)\phi_0 (x)dx
\]
for all $\psi\in C^\infty_0(\IR^d)$ and $\phi_1$ belonging to the domain of $\partial$. Therefore
\[
\int_\Omega\int_D A(x,\omega)(\nabla u_0(x)+U_1(x,\omega))\cdot\partial\phi_1(\omega)\psi(x) dx\IP(d\omega)=0.
\]
From this, we have for almost all $\omega\in\Omega$
\be
\int_D\int_{\IR^d}A(x,T(y)\omega)(\nabla u_0(x)+U_1(x,T(y)\omega))\cdot\nabla_y\phi_1(y))\psi(x)dydx=0
\label{eq:11}
\ee
$\forall\,\psi\in C^\infty_0(D), \phi_1\in C^\infty_0(\IR^d)$. On the other hand, as
\[
\int_\Omega\int_DA(x,\omega)(\nabla u_0(x)+U_1(x,\omega))\cdot\nabla\phi_0(x) dx\IP(d\omega)=\int_D f(x)\phi_0(x)dx,
\]
due to ergodicity
\[
\lim_{\rho\to \infty}{1\over |Y^\rho|}\int_{Y^\rho}\int_D A(x,T(y)\omega)(\nabla u_0(x)+U_1(x,T(y)\omega))\cdot\nabla\phi_0(x)dx dy=\int_D f(x)\phi_0(x)dx
\]
for almost all $\omega\in\Omega$. Thus almost surely $\forall\,\phi_1\in H^1_0(Y^\rho)$
\be
\lim_{\rho\to 0}{1\over|Y^\rho|}\int_D\int_{Y^\rho}A(x,T(y)\omega)(\nabla u_0(x)+U_1(x,T(y)\omega)\cdot(\nabla\phi_0(x)+\psi(x)\nabla_y\phi_1(y))dydx=\int_Df(x)\phi_0(x) dx.
\label{eq:u0U1}
\ee
As $U_1(x,T(y)\omega)$ is a potential vector with respect to $y$,
this leads us to consider the following approximation problem. 
We denote by $V_1^\rho:=L^2(D,H^1_0(Y^\rho))$ and $\bV^\rho=V\times V_1^\rho$
with the norm
\[
\|(w_0,w_1)\|_{\bV^\rho}=\|w_0\|_V+\|w_1\|_{V_1^\rho}
\]
for $(w_0,w_1)\in \bV^\rho$.
We define the bilinear form $B^\rho(\omega;(\cdot,\cdot),(\cdot,\cdot)):\bV^\rho\times\bV^\rho\to \IR$ by
\beqas
&&B^\rho(\omega;(w_0,w_1), (\phi_0,\phi_1))\\
&&={1\over |Y^\rho|}\int_D\int_{Y^\rho}A(x,T(y)\omega)(\nabla w_0(x)+\nabla_yw_1(x,y))\cdot(\nabla\phi_0(x)+\nabla_y\phi_1(x,y))dydx,
\eeqas
$\forall\,(w_0,w_1),(\phi_0, \phi_1)\in \bV^\rho$.
From \eqref{eq:boundednesscoerciveness}, it is straightforward to verify that $B^\rho(\omega;(\cdot,\cdot),(\cdot,\cdot))$ is  bounded and coercive for almost all $\omega\in\Omega$. We consider the approximate problem: Find $(u_0^\rho(\cdot,\omega), u_1^\rho(\cdot,\cdot,\omega))\in {\bV}^\rho$ such that
\be
B^\rho(\omega;(u_0^\rho(\cdot,\omega),u_1^\rho(\cdot,\cdot,\omega)),(\phi_0,\phi_1))=\int_Df(x)\phi_0(x)dx
\label{eq:approxprob}
\ee
$\forall\,(\phi_0, \phi_1)\in \bV^\rho$.  From the boundedness and the coercivity conditions of the bilinear form $B^\rho(\omega)$, this problem has a unique solution $(u_0^\rho,u_1^\rho)\in\bV^\rho$. 
We have the following estimate
\begin{lemma}
There is a positive constant $c$ which is independent of $\rho$ such that for almost all $\omega\in\Omega$
\[
\|u_0^\rho(\cdot,\omega)\|_V\le c,\ \ \|u_1^\rho(\cdot,\cdot,\omega)\|_{V_1^\rho}\le c|Y^\rho|^{1/2}.
\]
\end{lemma}
{\it Proof}\ \ Letting $\phi_0=u_0^\rho$ and $\phi_1=u_1^\rho$ in \eqref{eq:approxprob} and using \eqref{eq:boundednesscoerciveness}, we have
\[
c(\|u_0^\rho(\cdot,\omega)\|_V^2+|Y^\rho|^{-1}\|u_1^\rho(\cdot,\cdot,\omega)\|_{V_1^\rho}^2)\le \|f\|_{V'}\|u_0^\rho(\cdot,\omega)\|_V.
\]
We then get the conclusion.\hfill$\Box$

It is straighforward to verify that  $u_0^\rho$ satisfies the problem 
\be
-\nabla\cdot(A^{0\rho}(x,\omega)\nabla u_0^\rho(x,\omega))=f
\label{eq:u0rho}
\ee
with the Dirichlet boundary condition on $\partial D$, where $A^{0\rho}$ is defined in \eqref{eq:A0rho}. We then have the following approximation.
\begin{proposition} \label{prop:u0rhou0conv} Let $(u_0^\rho(\cdot,\omega), u_1^\rho(\cdot,\cdot,\omega))\in \bV^\rho$ be the solution of problem \eqref{eq:approxprob}. Then almost surely
\[
\lim_{\rho\to\infty}\|u_0-u_0^\rho(\cdot,\omega)\|_V=0.
\]
\end{proposition}
{\it Proof\ \ } 
First we show that  the coefficient $A^{0\rho}$ is uniformly bounded and coercive for almost all $(x,\omega)\in D\times\Omega$. As $|A_{ik}(x,\omega)|$ is uniformly bounded for all $x\in D$ and $\omega\in\Omega$, we have
\[
{1\over |Y^\rho|}\int_{Y^\rho} A_{ij}(x,T(y)\omega) dy
\]
is uniformly bounded for almost all $x\in D$ and $\omega\in\Omega$. From \eqref{eq:**}, we have
\beqas
{1\over |Y^\rho|}\int_{Y^\rho} A_{ik}(x,T(y)\omega){\partial N^{j\rho}\over\partial y_k}(x,y,\omega) dy
\le c{1\over |Y^\rho|}\int_{Y^\rho}\left|{\partial N^{j\rho}\over\partial y_k}(x,y,\omega)\right|dy\\
\le c{1\over |Y^\rho|}\left(\int_{Y^\rho}dy\right)^{1/2}\left(\int_{Y^\rho}\left|{\partial N^{j\rho}\over\partial y_k}(x,y,\omega)\right|^2dy\right)^{1/2}\le c.
\eeqas
For the coercivity, we note that
\beqas
A^{0\rho}_{ij}(x,\omega)={1\over |Y^\rho|}\int_{Y^\rho}A_{lk}(x,T(y)\omega)\left({\partial N^{j\rho}\over\partial y_k}+\delta_{kj}\right)\left({\partial N^{i\rho}\over\partial y_l}+\delta_{li}\right) dy
\eeqas
so for $\xi\in\IR^d$, 
\beqas
A_{ij}^{0\rho}\xi_i\xi_j\ge c{1\over|Y^\rho|}\sum_{k=1}^d\int_{Y^\rho}\left|\left({\partial N^{j\rho}\over\partial y_k}+\delta_{kj}\right)\xi_j\right|^2 dy\ge c|\xi|^2.
\eeqas
 Therefore $\|u_0^\rho(\cdot,\omega)\|_V$ is uniformly bounded.  From \eqref{eq:u0rho} and \eqref{eq:homprob}, we have
\[
-\nabla\cdot(A^0\nabla (u_0-u_0^\rho))=\nabla\cdot((A^0-A^{0\rho})\nabla u_0^{\rho}).
\]
Thus
\beqas
\int_DA^0(x)\nabla(u_0(x)-u_0^\rho(x,\omega))\cdot\nabla(u_0(x)-u_0^\rho(x,\omega))dx=\\
-\int_D(A^0-A^{0\rho})\nabla u_0^{\rho}(x,\omega)\cdot\nabla(u_0(x)-u_0^\rho(x,\omega))dx.
\eeqas
Therefore
\[
\|u_0-u_0^\rho(\cdot,\omega)\|_V\le c\|A^0-A^{0\rho}(\cdot,\omega)\|_{L^\infty(D)^{d\times d}}.
\]
We then get the conclusion.\hfill$\Box$

Approximation for $U_1(\cdot,T(\cdot)\omega)$ in terms of $u_1^\rho(\cdot,\cdot,\omega)$ is presented in Lemma \ref{lem:limrhoU1u1rho}.
%
\section{Finite element approximation of problem \eqref{eq:approxprob}}
\label{sec:FE}
We develop finite element approximations for the solution of problem \eqref{eq:approxprob} in this section. Due to the tensorized structure of the problem, we use tensor product FE approximations. We will first develop the full tensor product FE approximation which only  requires regularity for $u_1^\rho$ with respect to  the variables $x$ and $y$ separately, but the dimension of the FE space is exceedingly large. We then develop the sparse tensor product FE approximation which requires a stronger regularity condition to get the same level of accuracy but uses only  an essentially optimal number of degrees of freedom. We will show later that this stronger regularity condition is obtained under sufficient regularity conditions on the domain $D$, the coefficient $A$ and the forcing $f$.  Comparing to previous work on sparse tensor FEs for locally periodic problems such as \cite{HSelliptic} where the periodic unit cell has size 1, in this paper, we show that although the cube $Y^\rho$ has size $2\rho$, and the regularity norm of $u_1^\rho$ grows with $\rho$, the FE rate of convergence for $u_0^\rho$ is independent of $\rho$; it only depends on the mesh size. 

Let $D$ be a polygonal domain in $\IR^d$. We partition $D$ into hierarchical families of regular simplices. First the family ${\cal T}^1$ is obtained by dividing $D$ into triangular simplices of mesh size $O(2^{-1})$.  For $l=2,3,\ldots$, the family of simplices ${\cal T}^l$ is obtained by dividing each simplex in ${\cal T}^{l-1}$ into 4 congruent triangles when $d=2$ and 8 tetrahedra when $d=3$. The mesh size of simplices in ${\cal T}^l$ is $h_l=O(2^{-l})$. Similarly, we divide the cube $Y^\rho$ into families ${\cal T}^{\rho,l}$ of simplices of mesh size $h_l=O(2^{-l})$ (which does not depend on $\rho$). We define the following spaces:
\beqas
&&V^l=\{w\in H^1(D):\ w|_T\in P^1(K)\ \ \forall\,K\in{\cal T}^l\},\\
&&V^l_0=\{w\in H^1_0(D):\ w|_T\in P^1(K)\ \ \forall\,K\in{\cal T}^l\},\\
&&V^{\rho,l}_0=\{w\in H^1_0(Y^\rho):\ w|_T\in P^1(K)\ \ \forall\,K\in{\cal T}^{\rho,l}\}.
\eeqas
We then have the following approximation
\beqas
&&\inf_{w^l\in V^l}\|w-w^l\|_H\le ch_l|w|_{H^1(D)},\ \ \forall\,w\in H^1(D),\\
&&\inf_{w^l\in V^l_0}\|w-w^l\|_V\le ch_L|w|_{H^2(D)},\ \ \forall\,w\in H^2(D)\cap H^1_0(D),\\
&&\inf_{w^l\in V^{\rho,l}_0}\|w-w^l\|_{H^1_0(Y^\rho)}\le ch_L|w|_{H^2(Y^\rho)},\ \ \forall\,w\in H^2(Y^\rho)\cap H^1_0(Y^\rho),
\eeqas
where the constant $c$ in the last estimate does not depend on $\rho$ (see, e.g., \cite{Ciarlet} or \cite{BrennerScott}). 
As $u_1^\rho\in L^2(D, H^1_0(Y^\rho))\cong H\otimes H^1_0(Y^\rho)$, we employ tensor product FEs to approximate $u_1^\rho$. 
\subsection{Full tensor product FE approximation}
We consider the full tensor product FE space $V_1^{\rho,L}=V^L\otimes V_0^{\rho,L}$. We define the FE space $\bV^{\rho,L}=V^L\times V_1^{\rho,L}$. We consider the full tensor product FE problem: Find $(u_0^{\rho,L}(\cdot,\omega), u_1^{\rho,L}(\cdot,\cdot,\omega))\in\bV^{\rho,L}$ such that
\be
B^\rho(\omega;(u_0^{\rho,L}(\cdot,\omega),u_1^{\rho,L}(\cdot,\cdot,\omega)),(\phi_0^L,\phi_1^L))=\int_D f(x)\phi_0^L(x) dx,\ \ \forall\,(\phi_0^L,\phi_1^L)\in\bV^{\rho,L}.
\label{eq:fullprob}
\ee
We then have the following approximation.
\begin{lemma}\label{lem:Ceas}
There is a positive constant $c$ which is independent of $\rho$ such that for almost all $\omega\in \Omega$, the solution $(u_0^{\rho,L},u_1^{\rho,L})$ of the full tensor product approximating problem \eqref{eq:fullprob} satisfies
\beqas
|Y^\rho|\|u_0^\rho(\cdot,\omega)-u_0^{\rho,L}(\cdot,\omega)\|_V^2+\|u_1^\rho(\cdot,\cdot,\omega)-u_1^{\rho,L}(\cdot,\cdot,\omega)\|^2_{V_1^\rho}\\
\le c(|Y^\rho|\|u_0^\rho(\cdot,\omega)-\phi_0^L\|^2_V+\|u_1^\rho(\cdot,\cdot,\omega)-\phi_1^L\|_{V_1^\rho}^2)
\eeqas
$\forall\,(\phi_0^L,\phi_1^L)\in\bV^{\rho,L}$. 
\end{lemma}
{\it Proof}\ \ The proof of this lemma is straightforward. We have
\beqas
B^\rho(\omega;(u_0^\rho(\cdot,\omega)-u_0^{\rho,L}(\cdot,\omega),u_1^\rho(\cdot,\cdot,\omega)-u_1^{\rho,L}(\cdot,\cdot,\omega)),(u_0^\rho(\cdot,\omega)-u_0^{\rho,L}(\cdot,\omega),u_1^\rho(\cdot,\cdot,\omega)-u_1^{\rho,L}(\cdot,\cdot,\omega)) )\\
=B^\rho(\omega;(u_0^\rho(\cdot,\omega)-u_0^{\rho,L}(\cdot,\omega),u_1^\rho(\cdot,\cdot,\omega)-u_1^{\rho,L}(\cdot,\cdot,\omega)), (u_0^\rho(\cdot,\omega)-\phi_0^L,u_1^\rho(\cdot,\cdot,\omega)-\phi_1^L))
\eeqas
$\forall\,(\phi_0^L,\phi_1^L)\in\bV^{\rho,L}$. From \eqref{eq:boundednesscoerciveness}, we have 
\begin{align*}
&c_1\int_D\int_{Y^\rho}(|\nabla(u_0^\rho(x,\omega)-u_0^{\rho,L}(x,\omega))|^2+|\nabla_y(u_1^\rho(x,y,\omega)-u_1^{\rho,L}(x,y,\omega))|^2)dydx\\
&\le c_2\left(\int_D\int_{Y^\rho}(|\nabla(u_0^\rho(x,\omega)-u_0^{\rho,L}(x,\omega))|^2+|\nabla_y(u_1^\rho(x,y,\omega)-u_1^{\rho,L}(x,y,\omega)|^2)dydx\right)^{1/2}\\
&\cdot\left(\int_D\int_{Y^\rho}(|\nabla(u_0^\rho(x,\omega)-\phi_0^{L}(x))|^2+|\nabla_y(u_1^\rho(x,y,\omega)-\phi_1^{L}(x,y))|^2)dydx\right)^{1/2}.
\end{align*}
From this we get the conclusion.\hfill$\Box$

To quantify the error of the  approximate problem \eqref{eq:fullprob}, we define the following regularity space. Let ${\cal H}^\rho$ be the space $L^2(D, H^2(Y^\rho))\bigcap H^1(D, H^1_0(Y^\rho))$. We define the semi norm
\[
|w|_{{\cal H}^\rho}=\sum_{\alpha_1\in\IN_0^d, \atop |\alpha_1|=2}\left\|{\partial^{|\alpha_1|}w\over\partial y^{\alpha_1}}\right\|_{L^2(D\times Y^\rho)}+\sum_{\alpha_0,\alpha_1\in\IN_0^d\atop |\alpha_0|=1,|\alpha_1|\le 1}\left\|{\partial^{|\alpha_0|+|\alpha_1|}w\over\partial x^{\alpha_0}\partial y^{\alpha_1}}\right\|_{L^2(D\times Y^\rho)}.
\]
We then have the following approximating properties.
\begin{lemma}\label{lem:fullerror}
There is a positive constant $c$ independent of $\rho$ such that for all $w\in {\cal H}^\rho$
\[
\inf_{w^L\in V_1^{\rho,L}}\|w-w^L\|_{V_1^\rho}\le ch_L|w|_{{\cal H}^\rho}.
\]
\end{lemma}
The proof of this result is standard (see, e.g., \cite{BungartzGriebel, HSelliptic}). The fact that the constant $c$ in the lemma is independent of $\rho$ is because we can choose a constant $c$ independent of $\rho$ such that
\[
\inf_{v^L\in V_1^{\rho,L}}\|v^L-v\|_{H^1(Y^\rho)}\le ch_L|v|_{H^2(Y^\rho)}
\]
for all $v\in H^1_0(Y^\rho)\bigcap H^2(Y^\rho)$ (see, e.g.,\cite{Ciarlet}, \cite{BrennerScott}). 
\begin{lemma}\label{lem:fullreg}
Assume that $D$ is a convex domain or a domain of the $C^2$ class, $A\in C^1(\bar D, L^\infty(\Omega))$ and $\partial A\in C^1(\bar D, L^\infty(\Omega))^d$, and $f\in H$. Then $u_0^\rho\in H^2(D)$ and $u_1^\rho\in {\cal H}^\rho$. Further, there is a positive constant $c$ independent of $\rho$ such that
\[
\|u_0^\rho\|_{H^2(D)}\le c,\ \ |u_1^\rho|_{\cal H}\le c|Y^\rho|^{1/2}.
\]
\end{lemma}
We prove this lemma in Apprendix \ref{app:A}. From Lemmas \ref{lem:Ceas}, \ref{lem:fullerror} and \ref{lem:fullreg}, we have the following error estimate for the full tensor product FE approximating problem \eqref{eq:fullprob}.
\begin{proposition}\label{prop:fullerror}
Assume the hypothesis of Lemma \ref{lem:fullreg}. There is a positive constant $c$ independent of $\rho$ such that for almost all $\omega\in\Omega$, the solution of the approximating problem \eqref{eq:fullprob} satisfies
\[
\|u_0^\rho(\cdot,\omega)-u_0^{\rho,L}(\cdot,\omega)\|_V\le ch_L,\ \ \|u_1^\rho(\cdot,\cdot,\omega)-u_1^{\rho,L}(\cdot,\cdot,\omega)\|_{V_1}\le ch_L|Y^\rho|^{1/2}.
\]
\end{proposition}
\begin{remark} The dimension of the full tensor product FE space $\bV^{\rho,L}$ is $O(2^{2dL}\rho^d)$. 
\end{remark}
\subsection{Sparse tensor product FE approximation}
We develop in this section the sparse tensor product FE approximation which requires a far less number of degrees of freedom to achieve a prescribed level of accuracy. To this end, we define the orthogonal projection
\[
P^l:H\to V^l,\ \ P_1^{\rho,l}:H^1_0(Y^\rho)\to V_0^{\rho,l}
\]
with respect to the norm of $H$ and $H^1_0(Y^\rho)$ respectively. We define the detail spaces 
\[
W^l=(P^l-P^{l-1})V^{l},\ \ W^{\rho,l}=(P_1^{\rho,l}-P_1^{\rho,l-1})V^{\rho,l}_0,
\]
with the convention that $P^{-1}=0$ and $P^{\rho,-1}_1=0$. We note that
\[
V^L=\bigoplus_{l=1}^LW^l,\ \ V^{\rho,L}_0=\bigoplus_{l=1}^LW^{\rho,l}
\]
so
\[
V_1^{\rho,L}=\bigoplus_{l_1,l_2=1}^LW^{l_1}\otimes W^{\rho,l_2}.
\]
We define the sparse tensor product FE space as
\[
\hat V_1^{\rho,L}=\bigoplus_{l_1+l_2\le L}W^{l_1}\otimes W^{\rho,l_2}.
\]
Let $\hat{\bf V}^{\rho,L}=V^L\times\hat V_1^{\rho,L}$. We consider the sparse tensor product FE problem: Find $(\hat u_0^{\rho,L}(\cdot,\omega),\hat u_1^{\rho,L}(\cdot,\cdot,\omega))\in\hat{\bf V}^{\rho,L}$ such that
\be
B^\rho(\omega;(\hat u_0^{\rho,L}(\cdot,\omega),\hat u_1^{\rho,L}(\cdot,\cdot,\omega)),(\hat\phi_0^L,\hat\phi_1^L))=\int_Df(x)\hat\phi_0^L(x)dx,\ \ \forall\,(\hat\phi_0^L,\hat\phi_1^L)\in\hat{\bV}^{\rho,L}.
\label{eq:sparseprob}
\ee
An identical proof to that of Lemma \ref{lem:Ceas} shows that
\begin{lemma}\label{lem:sparseCeas}
There is a positive constant $c$ which is independent of $\rho$ such that for almost all $\omega\in \Omega$, the solution $(\hat u_0^{\rho,L},\hat u_1^{\rho,L})$ of the sparse tensor product FE approximating problem \eqref{eq:sparseprob} satisfies
\beqas
|Y^\rho|\|u_0^\rho(\cdot,\omega)-\hat u_0^{\rho,L}(\cdot,\omega)\|_V^2+\|u_1^\rho(\cdot,\cdot,\omega)-\hat u_1^{\rho,L}(\cdot,\cdot,\omega)\|^2_{V_1^\rho}\\
\le c(|Y^\rho|\|u_0^\rho(\cdot,\omega)-\hat\phi_0^L\|^2_V+\|u_1^\rho(\cdot,\cdot,\omega)-\hat\phi_1^L\|_{V_1^\rho}^2)
\eeqas
$\forall\,(\hat \phi_0^L,\hat\phi_1^L)\in\hat\bV^{\rho,L}$. 
\end{lemma}
We define the regularity space $\hat {\cal H}^\rho$ as the space of functions $w\in L^2(D,H^1_0(Y^\rho))$ such that for all $\alpha_0, \alpha_1\in\IN^d$ with $|\alpha_0|\le 1$ and $|\alpha_1|\le 2$, 
\[
{\partial^{|\alpha_0|+|\alpha_1|}w\over\partial x^{\alpha_0}\partial y^{\alpha_1}}\in L^2(D\times Y^\rho).
\]
In other words, $\hat{\cal H}^\rho=H^1(D, H^2(Y^\rho))$. We define the semi norm
\[
|w|_{\hat{\cal H}^\rho}=\sum_{\alpha_0,\alpha_1\in\IN^d\atop |\alpha_0|= 1,|\alpha_1|= 2}\left\|{\partial^{|\alpha_0|+|\alpha_1|}w\over\partial x^{\alpha_0}\partial y^{\alpha_1}}\right\|_{L^2(D\times Y^\rho)}.
\]
We have the following approximating result.
\begin{lemma}\label{lem:sparseerror}
There is a positive constant $c$ independent of $\rho$ such that for all $w\in\hat{\cal H}^\rho$
\[
\inf_{\hat w^L\in{\hat V_1^{\rho,L}}}\|w-\hat w^L\|_{V_1}\le cL^{1/2}h_L|w|_{\hat{\cal H}^\rho}.
\]
\end{lemma}
As for the proof of Lemma \ref{lem:fullerror}, the proof of this lemma is standard (see, e.g., \cite{BungartzGriebel, HSelliptic}).
We then have the following results on the regularity of the solution of problem \eqref{eq:approxprob}.
\begin{lemma}\label{lem:sparsereg}
Assume that $D$ is a convex domain or belongs to the $C^2$ class, $A\in C^1(\bar D, L^\infty(\Omega))$, $\partial A\in C^1(\bar D, L^\infty(\Omega))^d$, and $f\in H$. Then $u_0^\rho\in H^2(D)$ and $u_1^\rho\in \hat{\cal H}^\rho$. There is a positive constant $c$ independent of $\omega$ and $\rho$ such that
\[
\|u_0^\rho(\cdot,\omega)\|_{H^2(D)}\le c,\ \ |u_1^\rho(\cdot,\cdot,\omega)|_{\hat{\cal H}^\rho}\le c|Y^\rho|^{1/2}.
\]
\end{lemma}
We prove this lemma in Appendix \ref{app:A}. From Lemmas \ref{lem:sparseCeas}, \ref{lem:sparseerror} and \ref{lem:sparsereg}, we have the following error etimate.
\begin{proposition}\label{prop:sparseerror}
Assume the hypothesis of Lemma \ref{lem:sparsereg}, there is a positive constant $c$ independent of $\rho$ such that for almost all $\omega\in\Omega$, the solution of the sparse tensor product FE approximating problem \eqref{eq:sparseprob} satisfies
\[
\|u_0^\rho(\cdot,\omega)-\hat u_0^{\rho,L}(\cdot,\omega)\|_V\le cL^{1/2}h_L,\ \ \|u_1^\rho(\cdot,\cdot,\omega)-\hat u_1^{\rho,L}(\cdot,\cdot,\omega)\|_{V_1}\le cL^{1/2}h_L|Y^\rho|^{1/2}.
\]
\end{proposition}
\begin{remark}
The dimension of the sparse tensor product FE space $\hat {\bf V}^{\rho,L}$ is $O(L2^{dL}\rho^d)$ which is, apart from the multiplying factor $L$, equal to the dimension of the FE space for solving one truncated cell problem  in the cube $Y^\rho$. 
\end{remark}
\begin{remark}\label{rem:smoothdomains}
We note that tensor product FEs as considered above work equally in the case where the domain $D$ is not a polygon. For example, when $D$ is a smooth convex domain, to discretize $u_0^\rho$, we consider a polygon inscribed inside $D$ whose boundary edges (faces) are of the size $O(2^{-L})$. For constructing the sparse tensor product FEs to approximate $u_1$, we consider a polygon $\tilde D$ that contains $D$. The sparse tensor product FE spaces for approximating $u_1^\rho$ is the restriction of those for approximating functions in $L^2(\tilde D, H^1_0(Y^\rho))$ to $D\times Y^\rho$. We refer to \cite{Hoangmonotone} and \cite{XHelasticity} for details.

Alternatively, we can use parametric FEs as considered in \cite{HarbrechtSchwab}. 
\end{remark}
\begin{remark}
Brown and Hoang \cite{BrownHoang} develop an algorithm using multilevels of FE resolution and of numbers of Monte Carlo approximation samples at different macroscopic points to compute the expectation of the effective coefficient $A^{0\rho}(x,\omega)$ on the truncated domain $Y^\rho$ to approximate the homogenized coefficient $A^0(x)$. The algorithm needs an essentially optimal number of degrees of freedom for approximating the homogenized coefficients at a dense grid of macroscopic points $x$. 
\end{remark}


%
\section{Numerical corrector}
\label{sec:corrector}
We use the FE solutions of problem \eqref{eq:approxprob} developed in the previous sections to construct a numerical corrector for the solution of the stochastic two-scale problem \eqref{eq:2sprob}. 
To derive a numerical corrector for the solution of the two-scale problem \eqref{eq:2sprob}, we assume the following corrector result.
\be
\lim_{\ep\to 0}\IE[\|\nabla\ue-[\nabla u_0+w^i(\cdot,T({\cdot\over\ep})\omega)){\partial u_0\over\partial x_i}\|_{H^d}]=0.
\label{eq:corrector1}
\ee
%
\begin{remark}
Jikov et al. \cite{JKO} prove that when the coefficient $A$ does not depend on the slow variable, i.e. $A=A(\omega)$, then if $u_0\in C^\infty_0(D)$
\[
\lim_{\ep\to 0}\|\nabla\ue-[\nabla u_0+w^i(\cdot,T({\cdot\over\ep})\omega)){\partial u_0\over\partial x_i}\|_{H^d}=0,
\]
almost surely, which implies \eqref{eq:corrector1}. Examining the proof \cite{JKO}, we find that this result holds under the weaker condition that $u_0\in H^3_0(D)\cap W^{3,\infty}(D)$. We do not find a similar result in the literature for the case where the coefficient $A$ depends also on the slow variable $x$, i.e. $A=A(x,\omega)$. However, for quasi-periodic problems, assuming boundedness for the solution of the cell problems, we can derive this corrector result for the case where $A$ depends on the slow variable, following the standard procedure for the periodic case (\cite{BLP}, \cite{JKO}). 
\end{remark}
From \eqref{eq:2shomprob}, we have
\[
U_1(x,\omega)=w^i(x,\omega){\partial u_0\over\partial x_i}(x).
\]
Thus we can write \eqref{eq:corrector1} as
\be
\lim_{\ep\to 0}\IE[\|\nabla\ue-[\nabla u_0+U_1(\cdot,T({\cdot\over\ep})\omega)]\|_{H^d}]=0.
\label{eq:corrector}
\ee
 We use the FE approximations for $u_0^\rho$ and $u_1^\rho$ in the previous section to construct a numerical corrector in this section. We first note the following convergence result.
\begin{lemma}\label{lem:limrhoU1u1rho}
Almost surely, the solution of  problems \eqref{eq:2shomprob} and \eqref{eq:approxprob} satisfies
\[
\lim_{\rho\to \infty}{1\over|Y^\rho|^{1/2}}\|U_1(\cdot,T(\cdot)\omega)-\nabla_yu_1^\rho(\cdot,\cdot,\omega)\|_{L^2(D\times Y^\rho)}=0.
\]
\end{lemma}
{\it Proof}\ \ We consider the expression
\begin{align*}
&{\cal I}={1\over |Y^\rho|}\int_{Y^\rho}\int_D A(x,T(y)\omega)(\nabla(u_0(x)-u_0^\rho(x,\omega))+(U_1(x,T(y)\omega)-\nabla_yu_1^\rho(x,y,\omega)))\\
&\qquad\qquad\cdot(\nabla(u_0(x)-u_0^\rho(x,\omega))+(U_1(x,T(y)\omega)-\nabla_yu_1^\rho(x,y,\omega))) dx dy\\
&={1\over|Y^\rho|}\int_{Y^\rho}\int_D A(x,T(y)\omega)(\nabla u_0(x)+U_1(x,T(y)\omega))\cdot(\nabla u_0(x)+U_1(x,T(y)\omega))dxdy\\
&-{2\over|Y^\rho|}\int_{Y^\rho}\int_D A(x,T(y)\omega)(\nabla u_0(x)+U_1(x,T(y)\omega))\cdot(\nabla u_0^\rho(x,\omega)+\nabla_yu_1^\rho(x,y,\omega))dxdy\\
&+{1\over|Y^\rho|}\int_{Y^\rho}\int_DA(x,T(y)\omega)(\nabla u_0^\rho(x,\omega)+\nabla_yu_1^\rho(x,y,\omega))\cdot(\nabla u_0^\rho(x,\omega)+\nabla_yu_1^\rho(x,y,\omega))dxdy.
\end{align*}
From ergodicity property and \eqref{eq:2shomprob}, we have almost surely
\begin{align*}
&\lim_{\rho\to\infty}{1\over|Y^\rho|}\int_{Y^\rho}\int_DA(x,T(y)\omega)(\nabla u_0(x)+U_1(x,T(y)\omega))\cdot(\nabla u_0(x)+U_1(x,T(y)\omega))dx dy\\
&=\int_\Omega\int_D A(x,\omega)(\nabla u_0(x)+U_1(x,\omega))\cdot(\nabla u_0(x)+U_1(x,\omega))dx\IP(d\omega)\\
&=\int_D f(x)u_0(x)dx.
\end{align*}
Further
\begin{align*}
\lim_{\rho\to\infty}{1\over|Y^\rho|}\int_{Y^\rho}\int_D A(x,T(y)\omega)(\nabla u_0^\rho+\nabla_yu_1^\rho)\cdot(\nabla u_0^\rho+\nabla_yu_1^\rho)dxdy=\lim_{\rho\to\infty}\int_D fu_0^\rho dx=\int_D fu_0dx.
\end{align*}
We find the limit of 
\[
{1\over|Y^\rho|}\int_{Y^\rho}\int_D A(x,T(y)\omega)(\nabla u_0(x)+U_1(x,T(y)\omega))\cdot(\nabla u_0^\rho(x,\omega)+\nabla_yu_1^\rho(x,y,\omega))dxdy
\]
when $\rho\to \infty$. From \eqref{eq:11},
\[
{1\over |Y^\rho|}\int_{Y^\rho}\int_D A(x,T(y)\omega)(\nabla u_0(x)+U_1(x,T(y)\omega))\cdot\nabla_y u_1^\rho (x,y,\omega) dxdy= 0.
\]
We now show that
\[
\lim_{\rho\to\infty}{1\over |Y^\rho|}\int_{Y^\rho}\int_D A(x,T(y)\omega)(\nabla u_0(x)+U_1(x,T(y)\omega))\cdot\nabla u_0^\rho(x,\omega) dxdy=\int_D f(x)u_0(x) dx
\]
almost surely. From the ergodicity property and equation \eqref{eq:2shomprob}, we have
\beqas
\lim_{\rho\to\infty}  {1\over |Y^\rho|}\int_{Y^\rho}\int_D A(x,T(y)\omega)(\nabla u_0(x)+U_1(x,T(y)\omega))\cdot\nabla u_0(x) dxdy=\\
\int_\Omega\int_D A(x,\omega)(\nabla u_0(x)+U_1(x,\omega))\cdot\nabla u_0(x) dx \IP(d\omega)=\int_D f(x)u_0(x) dx.
\eeqas
We consider
\begin{align*}
&{1\over|Y^\rho|}\int_{Y^\rho}\int_D A(x,T(y)\omega)(\nabla u_0(x)+U_1(x,T(y)\omega))\cdot\nabla(u_0^\rho(x,\omega)-u_0(x))dxdy\\
&\le \left({1\over|Y^\rho|}\int_{Y^\rho}\int_D|\nabla u_0(x)+U_1(x,T(y)\omega)|^2 dxdy\right)^{1/2}\left({1\over |Y^\rho|}\int_{Y^\rho}\int_D|\nabla(u_0^\rho(x,\omega)-u_0(x))|^2 dxdy\right)^{1/2}.
\end{align*}
From \eqref{prop:u0rhou0conv},
\[
\lim_{\rho\to\infty}{1\over |Y^\rho|}\int_{Y^\rho}\int_D|\nabla(u_0^\rho(x,\omega)-u_0(x))|^2 dxdy=\lim_{\rho\to\infty}\int_D|\nabla(u_0^\rho(x,\omega)-u_0(x))|^2 dx=0.
\]
Due to ergodicity, we have
\begin{align*}
\lim_{\rho\to\infty}{1\over|Y^\rho|}\int_{Y^\rho}\int_D|\nabla u_0(x)+U_1(x,T(y)\omega)|^2 dxdy=\int_\Omega\int_D|\nabla u_0(x)+U_1(x,\omega)|^2 dxP(d\omega),
\end{align*}
so almost surely
\[
{1\over|Y^\rho|}\int_{Y^\rho}\int_D|\nabla u_0(x)+U_1(x,T(y)\omega)|^2 dxdy
\]
is bounded. Thus
\[
\lim_{\rho\to\infty} {1\over|Y^\rho|}\int_{Y^\rho}\int_D A(x,T(y)\omega)(\nabla u_0(x)+U_1(x,T(y)\omega))\cdot\nabla(u_0^\rho(x,\omega)-u_0(x))dxdy=0.
\]
Therefore
\[
\lim_{\rho\to\infty} {1\over|Y^\rho|}\int_{Y^\rho}\int_D A(x,T(y)\omega)(\nabla u_0(x)+U_1(x,T(y)\omega))\cdot\nabla u_0^\rho(x,\omega)dxdy=\int_D f(x)u_0(x)dx.
\]
Thus
\[
	\lim_{\rho\to\infty}{\cal I}=0
\]
almost surely. 
From \eqref{eq:boundednesscoerciveness}, we have 
\[
\lim_{\rho\to\infty}{1\over|Y^\rho|^{1/2}}\|U_1(\cdot,T(\cdot)\omega)-\nabla_yu_1^\rho(\cdot,\cdot,\omega)\|_{L^2(D\times Y^\rho)^d}=0.
\]
\hfill$\Box$

From now on, we denote the FE solution of both the full tensor product FE and sparse tensor product FE approximating problems \eqref{eq:fullprob} and \eqref{eq:sparseprob} by $u_0^{\rho,L}$ and $u_1^{\rho,L}$. 
We will use these together with \eqref{eq:corrector} to construct a numerical corrector for $\ue$. However, although $\nabla_yu_1^{\rho,L}(x,y,\omega)$ approximates $U_1(x,T(y)\omega)$ in $L^2(D\times Y^\rho)$, in general, $\nabla_yu_1^{\rho,L}(x,{x\over\ep},\omega)$ does not approximate $U_1(x,T({x\over\ep})\omega))$ in $H^d$. We thus define the following map ${\cal U}^\ep: L^1(D,L^1_{loc}(\IR^d))\to L^1(D)$ and use
it to construct the numerical corrector. For $i\in\IZ^d$, we denote the cube $i2\rho+Y^\rho$ by $Y^{i\rho}$ and the cube $\ep(i2\rho+Y^\rho)$ by $ Y^{i\ep\rho}$. Let $I\subset \IZ^d$ be the set of all $i\in \IZ^d$ such that $Y^{i\ep\rho}\cap D\ne\emptyset$. For $\Phi\in L^1(D, L^1_{loc} (\IR^d))$, we define
\[
{\cal U}^\ep(\Phi)(x)={1\over|Y^{i\ep\rho}|}\int_{Y^{i\ep\rho}}\Phi(z,{x\over\ep})dz,
\]
if $x\in Y^{i\ep\rho}$. The function $\Phi$ is understood as $0$ if $z\notin D$. We then have the following result.
\begin{lemma}
\label{lem:5.4}
For $\Phi\in L^1(D, L^1_{loc}(\IR^d))$, 
\[
\int_{D^\ep}{\cal U}^\ep(\Phi)(x)dx=\sum_{i\in I}{1\over(2\rho^d)}\int_{Y^{i\ep\rho}}\int_{Y^{i\rho}}\Phi(z,y)dydz
\]
where $D^\ep=\cup_{i\in I}Y^{i\ep\rho}$. 
\end{lemma}
{\it Proof}\ \ We have
\[
\int_{D^\ep}{\cal U}^\ep(\Phi)(x) dx=\sum_{i\in I}{1\over|Y^{i\ep\rho}|}\int_{Y^{i\ep\rho}}\int_{Y^{i\ep\rho}}\Phi(z,{x\over\ep}) dzdx.
\]
Making the change of variable $y=x/\ep$, we get the conclusion.\hfill$\Box$

For $x\in D$ and $y\in\IR^d$, we define the function $U_1^{\rho,L}(x,y,\omega)$ by
\[
U_1^{\rho,L}(x,y,\omega)=\nabla_y u_1^{\rho,L}(x,y-i2\rho,T(i2\rho)\omega)
\]
if $y\in Y^{i\rho}$. We have the following result.
\begin{lemma}\label{lem:EcalUU1U1rho}
For all $\ep>0$ we have
\begin{eqnarray*}
&&\IE\left[\int_D|{\cal U}^\ep(U_1(\cdot,T(\cdot)\omega)-U_1^{\rho,L}(\cdot,\cdot,\omega))(x)|^2 dx\right]\\
&&\qquad\qquad\qquad\le c{1\over|Y^\rho|}\IE\left[\int_D\int_{Y^\rho}|U_1(x,T(y)\omega)-\nabla_yu_1^{\rho,L}(x,y,\omega)|^2 dydx\right].
\end{eqnarray*}
\end{lemma}
{\it Proof}\ \ For all $\Phi\in L^2(D, L^2_{loc}(\IR^d))$, we have
\[
{\cal U}^\ep(\Phi)(x)^2={1\over|Y^{i\ep\rho}|^2}\left(\int_{Y^{i\ep\rho}}\Phi(z,{x\over\ep})dz\right)^2\le {1\over|Y^{i\ep\rho}|}\int_{Y^{i\ep\rho}}\Phi(z,{x\over\ep})^2dz={\cal U}^\ep(\Phi^2)(x).
\]
Thus
\begin{align*}
&{\cal J}:=\IE\left[\int_D|{\cal U}^\ep(U_1(\cdot,T(\cdot)\omega)-U_1^{\rho,L}(\cdot,\cdot,\omega))(x)|^2 dx\right]\\
&\le\IE\left[ \int_D{\cal U}^\ep((U_1(\cdot,T(\cdot)\omega)-U_1^{\rho,L}(\cdot,\cdot,\omega))^2)(x)dx\right]\\
&\le \IE\left[ \int_{D^\ep}{\cal U}^\ep((U_1(\cdot,T(\cdot)\omega)-U_1^{\rho,L}(\cdot,\cdot,\omega))^2)(x)dx\right].
\end{align*}
From Lemma \ref{lem:5.4}, using ergodicity
\begin{align*}
&{\cal J}\le\IE\left[{1\over(2\rho)^d}\sum_{i\in I}\int_{Y^{i\ep\rho}}\int_{Y^{i\rho}}((U_1(x,T(y)\omega)-U_1^{\rho,L}(x,y,\omega))^2dydx\right]\\
&=\sum_{i\in I}{1\over(2\rho)^d}\IE\left[\int_{Y^{i\ep\rho}}\int_{Y^{\rho}}((U_1(x,T(y)T(2\rho i)\omega)-\nabla_yu_1^{\rho,L}(x,y,T(2\rho i)\omega))^2 dydx\right]\\
&=\sum_{i\in I}{1\over(2\rho)^d}\IE\left[\int_{Y^{i\ep\rho}}\int_{Y^\rho}((U_1(x,T(y)\omega)-\nabla_yu_1^{\rho,L}(x,y,\omega))^2 dydx\right]\\
&={1\over(2\rho)^d}\IE\left[\int_D\int_{Y^\rho}((U_1(x,T(y)\omega)-\nabla_yu_1^{\rho,L}(x,y,\omega))^2 dydx\right]
\end{align*}
(note that $U_1(x,T(y)\omega)$ and $\nabla_yu_1^{\rho L}$ are understood as $0$ when $x\notin D$ in the definition of the operator ${\cal U}^\ep$. 
\hfill$\Box$

Next we show the following estimate
\begin{lemma}\label{lem:EU1calU1ep}
Assume that $u_0\in H^2(D)$. 
Then 
\[
\IE\left[\int_D|U_1(x,T({x\over\ep})\omega)-{\cal U}^\ep(U_1(\cdot,T(\cdot)\omega)(x)|^2 dx\right]\le c(\rho\ep).
\]
\end{lemma}
{\it Proof}\ \ 
Let $J$ be the set of the index $i\in I$ such that $Y^{i\ep\rho}\subset D$. We then have 
\begin{align*}
&\IE\left[\int_{\cup_{i\in J}Y^{i\ep\rho}}|U_1(x,T({x\over\ep})\omega)-{\cal U}^\ep(U_1(\cdot,T(\cdot)\omega)(x)|^2 dx\right]\\
&=\sum_{i\in J}\IE\left[\int_{Y^{i\ep\rho}}|U_1(x,T({x\over\ep})\omega)-{1\over |Y^{i\ep\rho}|}\int_{Y^{i\ep\rho}}U_1(z,T({x\over\ep})\omega)dz|^2dx\right]\\
&=\sum_{i\in J}\IE\left[\int_{Y^{i\ep\rho}}\left|{1\over |Y^{i\ep\rho}|}\int_{Y^{i\ep\rho}}(U_1(x,T({x\over\ep})\omega)-U_1(z,T({x\over\ep})\omega))dz\right|^2dx\right]\\
&\le \sum_{i\in J}{1\over |Y^{i\ep\rho}|}\int_{Y^{i\ep\rho}}\int_{Y^{i\ep\rho}}\IE\left[|U_1(x,T({x\over\ep})\omega)-U_1(z,T({x\over\ep})\omega))|^2\right] dzdx\\
&=\sum_{i\in J}{1\over|Y^{i\ep\rho}|}\int_{Y^{i\ep\rho}}\int_{Y^{i\ep\rho}}\IE[|U_1(x,\omega)-U_1(z,\omega)|^2] dz dx.
\end{align*}
We note that
\[
U_1(x,\omega)-U_1(z,\omega)=(w^i(x,\omega)-w^i(z,\omega)){\partial u_0\over\partial x_i}(x)+\left({\partial u_0\over\partial x_i}(x)-{\partial u_0\over\partial x_i}(z)\right)w^i(z,\omega).
\]
We have
\begin{align*}
&\int_{Y^{i\ep\rho}}\int_{Y^{i\ep\rho}}\left({\partial u_0\over\partial x_i}(x)\right)^2\IE[(w^i(x,\omega)-w^i(z,\omega))^2] dz dx\\
&\le c(\ep\rho)^2|Y^{i\ep\rho}|\int_{Y^{i\ep\rho}}\left({\partial u_0\over\partial x_i}(x)\right)^2dx
\end{align*}
due to \eqref{eq:Lipwi}.
Thus
\begin{align*}
&\sum_{i\in J}{1\over |Y^{i\ep\rho}|}\int_{Y^{i\ep\rho}}\int_{Y^{i\ep\rho}}\IE[(w^i(x,\omega)-w^i(z,\omega))^2]\left({\partial u_0\over\partial x_i}(x)\right)^2dz dx\le c(\ep\rho)^2.
\end{align*}
As $\IE[|w^i(z,\omega)|^2]$ is uniformly bounded with respect to $z$, we now show that
\[
\sum_{i\in J}{1\over |Y^{i\ep\rho}|}\int_{Y^{i\ep\rho}}\int_{Y^{i\ep\rho}}\left({\partial u_0\over\partial x_i}(x)-{\partial u_0\over\partial x_i}(z)\right)^2 dzdx\le c(\ep\rho)^2.
\]
As $u_0\in H^2(D)$, it is sufficient to to show that for $\phi \in H^1(D)$, 
\[
\sum_{i\in J}{1\over |Y^{i\ep\rho}|}\int_{Y^{i\ep\rho}}\int_{Y^{i\ep\rho}}\left(\phi(x)-\phi(z)\right)^2 dzdx\le c(\ep\rho)^2\|\phi\|_{H^1(D)}^2.
\]
The proof is similar to that of Lemma 5.5 in \cite{HSmultirandom}. We include the proof here for completeness. 
We have
\begin{align*}
&\sum_{i\in J}{1\over |Y^{i\ep\rho}|}\int_{Y^{i\ep\rho}}\int_{Y^{i\ep\rho}}\left(\phi(x)-\phi(z)\right)^2 dzdx\\
&\le c\sum_{i\in J}\sum_{j=1}^d{1\over |Y^{i\ep\rho}|}\int_{Y^{i\ep\rho}}\int_{Y^{i\ep\rho}}|\phi(x_1,\ldots,x_j,z_{j+1},\ldots,z_d)-\phi(x_1,\ldots,x_{j-1},z_j,\ldots,z_d)|^2 dzdx\\
&=c\sum_{i\in J}\sum_{j=1}^d{1\over|Y^{i\ep\rho}|}\int_{Y^{i\ep\rho}}\int_{Y^{i\ep\rho}}\left|\int_{z_j}^{x_j}{\partial\phi\over\partial t}(x_1,\ldots,x_{j-1},t,z_{j+1},\ldots,z_d)dt\right|^2 dzdx\\
&\le c(\ep\rho)\sum_{i\in J}\sum_{j=1}^d{1\over |Y^{i\ep\rho}|}\int_{Y^{i\ep\rho}}\int_{Y^{i\ep\rho}}\int_{\ep(i_j-2\rho)}^{\ep(i_j+2\rho)}\left({\partial\phi\over\partial t}(x_1,\ldots,x_{j-1},t,z_{j+1},\ldots,z_d)\right)^2 dtdzdx\\
&=c(\ep\rho)\sum_{i\in J}\sum_{j=1}^d{1\over |Y^{i\ep\rho}|}(\ep\rho)^{d+1}\left\|{\partial\phi\over\partial x_j}\right\|^2_{L^2(Y^{i\ep\rho})}\\
&\le c(\ep\rho)^2. 
\end{align*}
Let $D^{2\ep\rho}$ be the $2\ep\rho$ neighbourhood of $\partial D$. We have that $\cup_{i\in I\setminus J} Y^{i\ep\rho}\cap D\subset D^{2\ep\rho}$. Thus
\begin{align*}
&\IE\left[\int_{\cup_{i\in I\setminus J} Y^{i\ep\rho}\cap D}|U_1(x,T({x\over\ep})\omega)|^2\right]\le\int_{D^{2\ep\rho}}\IE[|w^i(x,\omega)|^2]\left({\partial u_0\over\partial x_i}(x)\right)^2 dx\\
&\le \sum_{i=1}^d\int_{D^{2\ep\rho}}\left({\partial u_0\over\partial x_i}(x)\right)^2 dx\le c\ep\rho
\end{align*}
as $u_0\in W^{1,\infty}(D)$. 
We further have
\begin{align*}
&\IE\left[\int_{\cup_{i\in I\setminus J}Y^{i\ep\rho}\cap D}|{\cal U}^\ep(U_1(\cdot,T(\cdot)\omega))(x)|^2dx\right]=\IE\left[\sum_{i\in I\setminus J}\int_{Y^{i\ep\rho}\cap D}{1\over |Y^{i\ep\rho}|^2}\left(\int_{Y^{i\ep\rho}}U_1(z,T({x\over\ep})\omega)dz\right)^2 dx\right]\\
&=\IE\left[\sum_{i\in I\setminus J}{1\over|Y^{i\ep\rho}|}\int_{Y^{i\ep\rho}\cap D}\int_{Y^{i\ep\rho}}U_1(z,T({x\over\ep})\omega)^2 dzdx\right]=\sum_{i\in I\setminus J}{|Y^{i\ep\rho}\cap D|\over|Y^{i\ep\rho}|}\int_{Y^{i\ep\rho}}\IE[U_1(z,\omega)^2]dz\\
&\le\sum_{i\in I\setminus J}\int_{Y^{i\ep\rho}\cap D}\IE[(w^i(x,\omega))^2]\left({\partial u_0\over\partial x_i}(z)\right)^2 dx\le c \int_{D^{2\ep\rho}}\left({\partial u_0\over\partial x_i}(z)\right)^2 dz\le c\ep\rho
\end{align*}
We then get the conclusion.\hfill$\Box$

We then have the following numerical corrector result.
\begin{theorem}\label{thm:numericalcorrector}
Assume 
the hypothesis of Lemmas \ref{lem:fullreg} and \ref{lem:sparsereg} hold, 
then for $\delta>0$, there are  constants $L_0>0$, $\rho_0>0$, $\ep_0>0$ such that if $L>L_0$, $\rho>\rho_0$, $\ep<\ep_0$ and $\ep<\delta/\rho$, then for the solution of the full and sparse tensor product FE approximating problems \eqref{eq:fullprob} and \eqref{eq:sparseprob}, we have
\[
\IE[\|\nabla\ue-[\nabla u_0^{\rho,L}+{\cal U}^\ep(U_1^{\rho,L})(\cdot)]\|_{H}]\le \delta.
\]
\end{theorem}
 %
 {\it Proof}\ \ We note that
 \begin{align*}
& \IE[\|\nabla\ue-[\nabla u_0^L+{\cal U}^\ep(U_1^{\rho,L})(\cdot)]\|_{H}]\le \IE[\|\nabla\ue-[\nabla u_0+U_1(\cdot, T({\cdot\over\ep})\omega)]\|_{H}]\\
 &+\IE[\|\nabla u_0-\nabla u_0^L\|_H]+\IE[\|U_1(\cdot,T({\cdot\over\ep})\omega)-{\cal U}^\ep(U_1(\cdot,T(\cdot)\omega))\|_H]\\
 &+\IE[\|{\cal U}^\ep(U_1(\cdot,T(\cdot)\omega))-{\cal U}^\ep(U_1^{\rho,L}(\cdot,\cdot,\omega))\|_H].
 \end{align*}
 From Lemma \ref{lem:EcalUU1U1rho}, we have
 \begin{align*}
 &\IE[\|{\cal U}^\ep(U_1(\cdot,T(\cdot)\omega))-{\cal U}^\ep(U_1^{\rho,L}(\cdot,\cdot,\omega)\|_H]^2\le {c\over|Y^\rho|}\IE[\|U_1(\cdot,T(\cdot)\omega)-\nabla_yu_1^\rho(\cdot,\cdot,\omega)\|_{L^2(D\times Y^\rho)}^2]\\
 &+{c\over|Y^\rho|}\IE[\|\nabla_yu_1^\rho (\cdot,\cdot,\cdot)-\nabla_yu_1^{\rho,L}(\cdot,\cdot,\cdot)\|_{L^2(D\times Y^\rho)}^2].
 \end{align*}
 From Proposition \ref{prop:fullerror}, Proposition \ref{prop:sparseerror}, equation \eqref{eq:corrector},  Lemma \ref{lem:limrhoU1u1rho}, Lemma \ref{lem:EU1calU1ep}, 
we get the conclusion.\hfill$\Box$
 \begin{remark} If the hypothesis of Lemmas \ref{lem:fullreg} and \ref{lem:sparsereg} does not hold, i.e. the regularity requirements on $u_1^\rho$ for the convergence rate of the full and sparse tensor product FE approximation do not hold, then we still have
 \[
 \lim_{L\to\infty} \|u_0^\rho(\cdot,\omega)-u_0^{\rho,L}(\cdot,\omega)\|_{V}+\|u_1^\rho(\cdot,\cdot,\omega)-u_1^{\rho,L}(\cdot,\cdot,\omega)\|_{V_1^\rho}=0
 \]
 albeit an explicit convergence rate in terms of the mesh size $h_L$ is not available (here we denote the numerical solution of both the full and sparse tensor product FE approximations as $u_0^{\rho,L}$ and $u_1^{\rho,L}$). The numerical corrector result in Theorem \ref{thm:numericalcorrector} still holds but the constant $L_0$ now depends on $\rho$. 
 \end{remark}
 \begin{remark} For random homogenization, there are  limited results on the homogenization rate of convergence (see, e.g., Yurinskii \cite{Yurinskii}, Gloria and Otto \cite{GloriaOtto}, Armstrong et al. \cite{Armstrongbook}). However, assume that we have a homogenization convergence rate
 \[
 \IE[\|\ue-u_0\|_V]
 \le c\ep^\gamma
 \]
 for $\gamma>0$, then with the convergence rate \eqref{eq:A0rhoconvergence} and Proposition \ref{prop:u0rhou0conv}, we have
 \[
 \IE[\|u^\ep-u_0^{\rho,L}\|_V]\le c(\rho^{-\beta/2}+h_L+\ep^\gamma)
 \]
 for the solution of the full tensor product FE approximation problem \eqref{eq:fullprob}; and
 \[
 \IE[\|u^\ep-\hat u_0^{\rho,L}\|_V]\le c(\rho^{-\beta/2}+L^{1/2}h_L+\ep^\gamma)
 \]
 for the solution of the sparse tensor product FE approximation problem \eqref{eq:sparseprob}. 
 \end{remark}
\section{Numerical experiments}
\label{sec:numerical}
We first consider a quasi-periodic problem on the one dimensional domain $D=(0,1)$. Let the probability space $\Omega$ be the square $(0,1)\times (0,1)$ in $\IR^2$ with the Lebesgue probability measure. Let $\lambda$ be the vector $(1, \sqrt{2})^\top$. The dynamical system $T:\IR\to\Omega$ is defined as 
\[
T(y)\omega=(\omega+\lambda y)\, {\rm mod}\, 1=\left(\begin{matrix} (\omega_1+y)\,{\rm mod}\,1\\(\omega_2+\sqrt{2}y)\,{\rm mod}\,1\end{matrix}\right).
\]
We consider the coefficient 
\[
A(x,\omega)=(1+x)(3+\cos(2\pi\omega_1)+\cos(2\pi\omega_2))
\]
for $\omega=(\omega_1,\omega_2)\in\Omega$. We note that for a differentiable function $\phi$ belonging to the domain of $\partial_\omega$, 
\[
\partial_\omega\phi(\omega)=\lim_{y\to 0}{\phi(\omega_1+y,\omega_2+\sqrt{2}y)-\phi(\omega_1,\omega_2)\over y}={\partial\phi\over\partial\omega_1}(\omega)+\sqrt{2}{\partial\phi\over\partial\omega_2}(\omega).
\]
The cell problem \eqref{eq:abstractcell} is of the form
\[
\int_\Omega A(x,\omega)\left(1+{\partial N\over\partial\omega_1}+\sqrt{2}{\partial N\over\partial\omega_2}\right)\left({\partial\phi\over\partial\omega_1}+\sqrt{2}{\partial\phi\over\partial\omega_2}\right)d\omega=0,\ \ 
\forall\,\phi\in C^1_{per}(\Omega).
\]
We solve this cell problem numerically by FEs with a high accuracy level to compute the homogenized coefficient $A^0(x)$ from \eqref{eq:A0}. We find that approximately
\[
A^0(x)=2.6085(1+x).
\]
We choose the function $f$ in \eqref{eq:2sprob} such that the homogenized problem \eqref{eq:homprob} has the exact solution $u_0=x-x^2$. Problem \eqref{eq:approxprob} is solved by sparse tensor product FEs. We record in Table \ref{t:1} the numerical errors for $\IE[\|u_0^{\rho,L}-u_0\|_V^2]$ computed by using a highly accurate Gauss-Legendre quadrature rule in $\Omega$. 
\begin{table}[h!]
\centering
\begin{tabular}{ | c | c |c|c|c| } 
 \hline
 mesh level $L$ & error  for $\rho=1$ & error for $\rho=2$ & error for $\rho=3$ & error for $\rho=4$ \\ 
 \hline
 1&   0.08643387 & 0.08642242 & 0.08640925 & 0.08639875 \\
 2& 0.02175469 & 0.02168175 & 0.02163293 &  0.02160448\\
 3&  0.00547390 & 0.00541613 & 0.00535287 & 0.00530344\\
 4& 0.00153847 & 0.00147204 & 0.00139897 & 0.00134186\\
 5&  0.00056918 & 0.00050014 & 0.00042423 & 0.00036492 \\
 6& 0.00032794 & 0.00025823 & 0.00018158 & 0.00012169\\
 7& 0.00026771 & 0.00019782 & 0.00012098 & 0.00006095 \\
 8&  0.00025265 & 0.00018272 & 0.00010584 & 0.00004577\\
  9& 0.00024889 & 0.00017895 & 0.00010205 & 0.00004198\\
 \hline
\end{tabular}
\caption{Sparse tensor product FE error $\IE[\|u_0^{\rho,L}-u_0\|_V^2]$ for the one dimensional problem}
\label{t:1}
\end{table}
The numerical results show that for a fixed $\rho$, when we reduce the mesh size, the error reduces as the theoretical error estimate for the sparse tensor product FE approximations at first until the error due to the truncated cube $Y^\rho$ becomes more significant than the FE error. However, when we increase $\rho$, the error for smaller mesh sizes, where the FE error is less significant, decreases. This shows that the approximation gets more accurate when we enlarge the cube $Y^\rho$. 

We then consider a two-scale problem in the two dimensional domain $D=(0,1)\times (0,1)$.  Let $\mu$ be the matrix
\[
\mu=\left(\begin{matrix}
 1 & \sqrt{2} \\
 \sqrt{2} & 1
 \end{matrix}\right).
 \]
We consider the dynamical system $T:\IR^2\to \Omega$ such as
\[
T(y)\omega=(\omega+\mu y)\,{\rm mod}\, 1=\left(\begin{matrix} (\omega_1+y_1+\sqrt{2}y_2)\,{\rm mod}\,1\\ (\omega_2+\sqrt{2}y_1+y_2)\,{\rm mod}\,1\end{matrix}\right).
\]
The coefficient $A(x,\omega)$ for $x=(x_1,x_2)\in D$ and $\omega\in \Omega$ is
\[
A(x,\omega)=(1+x_1)(1+x_2)(3+\cos(2\pi\omega_1)+\cos(2\pi\omega_2)).
\]
For a function $\phi$ belonging to the domain of $\partial$, we have
\[
\partial_1\phi ={\partial\phi\over\partial\omega_1}+\sqrt{2}{\partial\phi\over\partial\omega_2},
\ \ \mbox{and} 
\ \ 
\partial_2\phi =\sqrt{2}{\partial\phi\over\partial\omega_1}+{\partial\phi\over\partial\omega_2}.
\]
For the cell problems \eqref{eq:abstractcell}, we solve the problems
\[
\int_\Omega A\left[\left(1+{\partial N^1\over\partial\omega_1}+\sqrt{2}{\partial N^1\over\partial\omega_2}\right)\left({\partial\phi\over\partial\omega_1}+\sqrt{2}{\partial\phi\over\partial\omega_2}\right)+\left(\sqrt{2}{\partial N^1\over\partial\omega_1}+{\partial N^1\over\partial\omega_2}\right)\left(\sqrt{2}{\partial\phi\over\partial\omega_1}+{\partial\phi\over\partial\omega_2}\right)\right]d\omega=0
\]
$\forall\,\phi\in C^1_{per}(\Omega)$,
and 
\[
\int_\Omega A\left[\left({\partial N^2\over\partial\omega_1}+\sqrt{2}{\partial N^2\over\partial\omega_2}\right)\left({\partial\phi\over\partial\omega_1}+\sqrt{2}{\partial\phi\over\partial\omega_2}\right)+\left(1+\sqrt{2}{\partial N^2\over\partial\omega_1}+{\partial N^2\over\partial\omega_2}\right)\left(\sqrt{2}{\partial\phi\over\partial\omega_1}+{\partial\phi\over\partial\omega_2}\right)\right]d\omega=0
\]
$\forall\,\phi\in C^1_{per}(\Omega)$,
numerically by FEs. We find that the homogenized coefficient $A^0(x)$ is approximately
\[
A^0(x)=(1+x_1)(1+x_2)\left(\begin{matrix}2.82491 & -0.16755\\-0.16755 & 2.82491\end{matrix}\right).
\]
We then choose the function $f$ so that the homogenized equation \eqref{eq:homprob} has the exact solution 
\[
u_0(x)=(x_1-x_1^2)(x_2-x_2^2)
\]
for $x=(x_1,x_2)\in D$.  We have the following numerical result.
 \begin{table}[h!]
\centering
\begin{tabular}{ | c | c |c|c| } 
 \hline
 mesh level $L$ & error  for $\rho=1$ & error for $\rho=2$ & error for $\rho=3$ \\ 
 \hline
 1&   0.00600307 & 0.00600258 & 0.00600243 \\
 2&   0.00143868 & 0.00143325 & 0.00143156 \\
 3&   0.00035356 & 0.00035110 & 0.00035004 \\
 4&   0.00008871 & 0.00008734 & 0.00008708 \\
 5&   0.00002200 & 0.00002191 & 0.00002175 \\
 \hline
\end{tabular}
\caption{Sparse tensor product FE error $\IE[\|u_0^{\rho,L}-u_0\|_V^2]$ for the two dimensional problem}
\label{t:2}
\end{table}
In the chosen range of the mesh sizes we experiment, the sparse tensor product FE approximation on the truncated domain converges to the exact solution $u_0$. The FE error dominates the error due to the truncated cube $Y^\rho$ so that the total error behaves as the theoretical error proved above for the sparse tensor product FE approximations. For the same FE mesh size, the total error is lightly reduced when the size of the cube $Y^\rho$ increases. 

The sparse tensor product FE convergence rate established in the paper holds when the coefficient $A$ is sufficiently smooth so that the solution of the truncated stochastic two-scale homogenized equation is sufficiently regular. However, when $A$ is not continuous, the sparse tensor product FEs may still work when the FE mesh fits into the surface of discontinuity. Now we consider the checker board problem where the two-scale coefficient $A$ is discontinuous. The problem is studied in details in \cite{JKO} where the homogenized coefficient can be computed explicitly in two dimensions. The two dimensional space $\IR^2$ is split into squares of size 1 whose centres are of integer components. Let $\cal M$ be the probability space of functions that takes values $1$ or $2$ in each square with probability $0.5$. Let $\lambda$ be the probability measure on $\cal M$. This probability space is ergodic with respect to the integer shift. We then define the probability space $\Omega$ to be the one of all piecewise constant functions obtained from those in $\cal M$ by a shift with a vector belonging to the unit cube $Y^1$, i.e. 
\[
\Omega=\{u(t):\ u(t)=\phi(t+\eta),\ \phi(t)\in{\cal M},\ \eta\in Y^1\}.
\]
The probability space $\Omega$ is associated with the product space ${\cal M}\times Y^1$ with the product probability measure. It is invariant and ergodic with respect to the shift operator. For $\omega=u(t)\in\Omega$, we define the coefficient
\[
A(\omega)=\left\{
\begin{array}{rl}
1 & \text{if } u(0)=1,\\
2 & \text{if } u(0)=2.
\end{array} \right.
\]
The coefficient $A$ is of the checker board type with value $1$ or $2$ in each square. The homogenized coefficient is $A^0(x)=2^{1/2}$.  We choose $f$ so that the homogenized equation has exact solution $u_0(x)=(x_1-x_1^2)(x_2-x_2^2)$ for $x=(x_1,x_2)\in D=(0,1)\times (0,1)$. The average over the probability space $\Omega$ involving taking the integral with respect to the shift $\eta\in Y^1$ and average over all the realizations of $\phi$ in $\cal M$. We solve equation \eqref{eq:approxprob} with $\rho=1, 2$ and $3$ by the sparse tensor product FEs. Integral over the unit cube $Y^1$ is  approximated by the Newton-Cotes rule with 5 quadrature points on each direction. When $\rho=1$, average over $\cal M$ is computed exactly. When $\rho=2$ and $3$, average over $\cal M$ is approximated with the Monte Carlo procedure with 1000 samples. Table \ref{t:3} presents the error for  $\IE[\|u_0^{\rho,L}-u_0\|_V^2]$. The results clearly indicates that the approximation gets better with larger values of $\rho$. The reduction rate of the error gets worse with smaller mesh sizes; this shows the more prominent effect of the domain truncating error when the FE mesh size gets smaller.
 \begin{table}[h!]
\centering
\begin{tabular}{ | c | c |c|c| } 
 \hline
 mesh level $L$ & error  for $\rho=1$ & error for $\rho=2$ & error for $\rho=3$ \\ 
 \hline
 2&   0.0632550 & 0.003578 & 0.003514 \\
 3&   0.0384572 & 0.001048 & 0.000974 \\
 4&   0.0274870 & 0.000372 & 0.000296 \\
 \hline
\end{tabular}
\caption{Sparse tensor product FE error $\IE[\|u_0^{\rho,L}-u_0\|_V^2]$ for the two dimensional checker board problem}
\label{t:3}
\end{table}
\begin{appendices}
\section{}\label{app:B}
We show Proposition \ref{prop:A0conv}. Let $N^{i\rho}(x,y,\omega)$ be the solution of the cell problem
\be
-\nabla_y\cdot(A(x,T(y)\omega)(\nabla_yN^{i\rho}(x,y,\omega)+e^i))=0
\label{eq:celli}
\ee
with the Dirichlet boundary condition $N^{i\rho}(x,y,\omega)=0$ when $y\in \partial Y^\rho$; $e^i$ is the $i$th unit vector in the standard basis of $\IR^d$. Following \cite{BourgeatPiatnitski}, we let $M^i(x,z,\omega)={1\over\rho}N^{i\rho}(x,\rho z,\omega)$. Then $M^i(x,z,\omega)$ is the solution of the problem 
\[
-\nabla_z\cdot(A(x,T(\rho z)\omega)(\nabla_zM^{i\rho}(x,z,\omega)+e^i))=0,\ \ \mbox{in}\ Y^1
\]
with the Dirichlet boundary condition on $\partial Y^1$. 
By the homogenization theory, almost surely, $M^{i\rho}\wc M^{i0}$ in $H^1_0(Y^1)$ when $\rho\to\infty$ where $M^{i0}$ is the solution of the homogenized equation
\be
-\nabla_z\cdot(A^0(x)(\nabla_z M^{i0}(x,z)+e^i))=0,\ \ \mbox{in}\ Y^1 
\label{eq:Mi0}
\ee
with the Dirichlet boundary condition on $\partial Y^1$; $A^0$ is the homogenized coefficient in \eqref{eq:A0}. Almost surely, we have 
\[
\int_{Y^1}A(x,T(\rho z)\omega)(\nabla_zM^{i\rho}(x,z,\omega)+e^i)dz\to \int_{Y^1}A^0(x)(\nabla_z M^{i0}(x)+e^i)dz
\]
when $\rho\to\infty$.
As the solution $M^{i0}=0$, and 
\[
\int_{Y^1}A(x,T(\rho z)\omega)(\nabla_zM^{i\rho}(x,z,\omega)+e^i)dz={1\over\rho^d}\int_{Y^\rho}A(x,T(y)\omega)(\nabla_yN^{i\rho}(x,y,\omega)+e^i)dy,
\]
we have that
\[
{1\over|Y^\rho|}\int_{Y^\rho}A(x,T(y)\omega)(\nabla_yN^{i\rho}(x,y,\omega)+e^i)dy\to{1\over |Y^1|}\int_{Y^1}A^0(x)(\nabla_z M^{i0}(x,z)+e^i)dz= A^0(x)e^i,
\]
when $\rho\to\infty$. We now show that almost surely the convergence is uniform with respect to $x$. First we show that $A^0(x)$ is uniformly continuous in $\bar D$. From \eqref{eq:abstractcell}, for $x, x'\in\bar D$, we have
\[
\int_\Omega A(x,\omega)(w^i(x,\omega)+e^i)\cdot(w^i(x,\omega)-w^i(x',\omega))\IP(d\omega)=0
\]
and 
\[
\int_\Omega A(x',\omega)(w^i(x',\omega)+e^i)\cdot(w^i(x,\omega)-w^i(x',\omega))\IP(d\omega)=0.
\]
Thus
\beqas
&&\int_\Omega A(x,\omega)(w^i(x,\omega)-w^i(x',\omega))\cdot(w^i(x,\omega)-w^i(x',\omega))\IP(d\omega)\\
&&=-\int_\Omega [(A(x,\omega)-A(x',\omega))w^i(x',\omega)\cdot(w^i(x,\omega)-w^i(x',\omega))\\
&&-(A(x,\omega)-A(x',\omega))e^i\cdot(w^i(x,\omega)-w^i(x',\omega))]\IP(d\omega).
\eeqas
From \eqref{eq:abstractcell}, we have that $w^i(x,\cdot)$ is uniformly bounded in $L^2(\Omega)$ with respect to $x$. Thus
\be
\|w^i(x,\cdot)-w^i(x',\cdot)\|_{L^2(\Omega)}\le c\|A(x,\cdot)-A(x',\cdot)\|_{L^\infty(\Omega)}\le c|x-x'|.
\label{eq:Lipwi}
\ee
From \eqref{eq:A0}
\beqas
A^0_{ij}(x)-A^0_{ij}(x')=\int_\Omega[(A_{ik}(x,\omega)-A_{ik}(x',\omega))\delta_{kj}+A_{ik}(x,\omega)(w^j_k(x,\omega)-w^j_k(x',\omega))\\
+(A_{ik}(x,\omega)-A_{ik}(x',\omega))w^j_k(x',\omega)]\IP(d\omega).
\eeqas
Thus
\[
|A^0_{ij}(x)-A^0_{ij}(x')|\le c|x-x'|.
\]
Next we show that $A^{0\rho}(\cdot,\omega)$ is also uniformly Lipschitz with respect to $\rho$ and $\omega$. From \eqref{eq:cellprob}, we have
\[
-\nabla_y\cdot(A(x,T(y)\omega)\nabla_yN^{i\rho}(x,y,\omega))=\nabla_y\cdot(A(x,T(y)\omega)e^i),
\]
so
\be
\|\nabla_yN^{i\rho}(x,\cdot,\omega)\|_{L^2(Y^\rho)}\le c\|A(x,T(\cdot)\omega)\|_{L^2(Y^\rho)}\le c|Y^\rho|^{1/2}.
\label{eq:**}
\ee
 From \eqref{eq:cellprob}, for $x,x'\in \bar D$ we have
\begin{eqnarray}
&&-\nabla_y\cdot(A(x,T(y)\omega)\nabla_y(N^{i\rho}(x,y,\omega)-N^{i\rho}(x',y,\omega)))
=\nabla_y\cdot((A(x,T(y)\omega)-A(x',T(y)\omega))e^i)\nonumber\\
&&+\nabla_y\cdot((A(x,T(y)\omega)-A(x',T(y)\omega))\nabla_yN^{i\rho}(x',T(y)\omega)).
\label{eq:***}
\end{eqnarray}
Therefore
\begin{eqnarray}
\|\nabla_y(N^{i\rho}(x,\cdot,\omega)-N^{i\rho}(x',\cdot,\omega))\|_{L^2(Y^\rho)^d}\le c\|A(x,T(\cdot)\omega)-A(x',T(\cdot)\omega)\|_{L^2(Y^\rho)}\nonumber\\
+\|A(x,T(\cdot)\omega)-A(x',T(\cdot)\omega)\|_{L^\infty(Y^\rho)}\|\nabla_yN^\rho(x,T(\cdot)\omega)\|_{L^2(Y^\rho)}
\le c|x-x'||Y^\rho|^{1/2}.
\label{eq:NrhoLipschitz}
\end{eqnarray}
From \eqref{eq:A0rho}, we have
\beqas
&&A^{0\rho}_{ij}(x,\omega)-A^{0\rho}_{ij}(x',\omega)={1\over|Y^\rho|}\int_{Y^\rho}(A_{ik}(x,T(y)\omega)-A_{ik}(x',T(y)\omega))\delta_{kj}dy\\
&&\qquad\qquad+{1\over|Y^\rho|}\int_{Y^\rho}[A_{ik}(x,T(y)\omega){\partial N^{j\rho}\over\partial y_k}(x,T(y)\omega)-A_{ik}(x',T(y)\omega){\partial N^{j\rho}\over\partial y_k}(x',T(y)\omega)]dy\\
&&\qquad\qquad={1\over|Y^\rho|}\int_{Y^\rho}[A_{ik}(x,T(y)\omega)-A_{ik}(x',T(y)\omega)]\delta_{kj}dy\\
&&\qquad\qquad+{1\over|Y^\rho|}\int_{Y^\rho}\Big[A_{ik}(x,T(y)\omega)-A_{ik}(x',T(y)\omega)){\partial N^{j\rho}\over\partial y_k}(x,T(y)\omega)\\
&&\qquad\qquad+A_{ik}(x',T(y)\omega)({\partial N^{j\rho}\over\partial y_k}(x,T(y)\omega)-{\partial N^{j\rho}\over\partial y_k}(x',T(y)\omega))\Big]dy.
\eeqas
From \eqref{eq:NrhoLipschitz} and the fact that $A(\cdot,\omega)$ is uniformly bounded in $C^{0,1}(\bar D)$, we have that $A^{0\rho}_{ij}(\cdot,\omega)$ is uniformly Lipschitz in $\bar D$.
For each $n\in\IN$, we consider the cubes $D_i^n$ ($i=1,\ldots,M_n$) of size $1/n$ with vertices having coordinates of the form $j/n$ for $j\in\IZ^d$ that intersects $D$. In each cube $D_i^n$ for $i=1,\ldots,M_n$  we choose a point $x_i^n\in D$. Then there is a set $\Omega_n\subset\Omega$ with $\IP(\Omega_n)=1$ such that for all $\omega\in\Omega_n$ $A^{0\rho}(x_i^n,\omega)\to A^0(x_i^n)$ when $\rho\to\infty$ for all $i=1,\ldots,n$. For each $x\in D_i^n$, $|A^0(x)-A^0(x_i^n)|\le c/n$ and $|A^{0\rho}(x,\omega)-A^{0\rho}(x_i^n,\omega)|\le c/n$. For each $\omega\in\Omega_n$, there is a constant $\rho_n=\rho_n(\omega)$ such that when $\rho>\rho_n$
\[
|A^{0\rho}(x_i^n,\omega)-A^0(x_i^n)|\le 1/n
\]
for all $i=1,\ldots,M_n$. Thus for $\omega\in \Omega_n$, for all $x\in D_i^n$ for all $i$, when $\rho>\rho_n$
\beqas
|A^{0\rho}(x,\omega)-A^0(x)|\le |A^{0\rho}(x,\omega)-A^{0\rho}(x_i^n,\omega)|+|A^{0\rho}(x_i^n,\omega)-A^{0}(x_i^n)|+|A^{0}(x_i^n)-A^0(x)|\le  c/n.
\eeqas
Let $\Omega_0=\cap_{n=1}^\infty\Omega_n$. For $\omega\in\Omega_0$, $\lim_{\rho\to\infty}\|A^{0\rho}(\cdot,\omega)-A^0(\cdot)\|_{L^\infty(D)}=0$.

\section{}\label{app:A}
{\it Proof of Lemma \ref{lem:sparsereg}}.  We show Lemma \ref{lem:sparsereg} (which implies Lemma \ref{lem:fullreg}). 
From \eqref{eq:u0rho}, we deduce that 
\[
\|u_0^\rho(\cdot,\omega)\|_{H^2(D)}\le c\|f\|_{H}
\]
where the constant $c$ only depends on the convex domain $D$ and the $C^{0,1}(\bar D)$ norm of $A^{0\rho}(\cdot,\omega)$ (see Grisvard \cite{Grisvard} Theorems 3.1.3.1 and 3.2.1.2). 

We note that
\[
u_1^\rho(x,y,\omega)={\partial u_0^\rho\over\partial x_i}(x,\omega)N^i(x,y,\omega)
\]
where $N^i(x,y,\omega)$ is the solution of the cell problem \eqref{eq:celli}.
To show that $|u_1^\rho|_{\hat{\cal H}}\le c|Y^\rho|^{1/2}$, it is sufficient to show that 
\[
|N^{i\rho}(x,\cdot,\omega)|_{H^2(Y))}\le c|Y^\rho|^{1/2},\ \mbox{and}\ |N^{i\rho}(x,\cdot,\omega)-N^{i\rho}(x',\cdot,\omega)|_{H^2(Y^\rho)}\le c|x-x'||Y^\rho|^{1/2}
\]
for all $x, x'\in D$.
From \eqref{eq:celli}
\[
-\nabla_y\cdot(A(x,T(y)\omega)\nabla_yN^{i\rho}(x,y,\omega))=\nabla_y\cdot(A(x,T(y)\omega)e^i).
\]
From the proof of Lemma 3.1.3.2 of \cite{Grisvard}, there is a constant $c$ that only depends on the Lipschitz norm of the coefficient $A(x,T(\cdot)\omega)$ so that the following inequality on the semi $H^2(Y^\rho)$ norm holds:
\[
|N^{i\rho}(x,\cdot,\omega|_{H^2(D)}\le c \|\nabla_y\cdot(A(x,T(\cdot)\omega)e^i)\|_{L^2(Y^\rho)}+\|N^{i\rho}(x,\cdot,\omega)\|_{H^1(Y^\rho)}).
\]
Indeed, Grisvard \cite{Grisvard} proves this estimate for the full $H^2$ norm so the constant depends on the diameter of the domain. However, from the proof of Theorem 3.1.2.1 of \cite{Grisvard}, we have that for the semi $H^2$ norm only, this constant does not depend on the diameter of the domain (the only part of the proof in \cite{Grisvard} where this constant depends on the domain is when using the Poincare inequality to bound the $L^2$ norm of the solution). Thus using \eqref{eq:**}, we have
\be
|N^{i\rho}(x,\cdot,\omega)|_{H^2(Y^\rho)}\le c|Y^\rho|^{1/2}.
\label{eq:Nirho}
\ee
For $x, x'\in\bar D$ from \eqref{eq:***}, we have
\beqas
&&\nabla_y\cdot(A(x,T(y)\omega)\cdot\nabla_y(N^{i\rho}(x,y,\omega)-N^{i\rho}(x',y,\omega)))=\\
&&\nabla_y\cdot((A(x,T(y)\omega)-A(x',T(y)\omega))e^i)-\nabla_y\cdot((A(x,T(y)\omega)-A(x',T(y)\omega))\nabla_yN^{i\rho}(x',T(y)\omega)).
\eeqas
From \eqref{eq:NrhoLipschitz} and \eqref{eq:Nirho}, we have
\[
|N^{i\rho}(x,\cdot,\omega)-N^{i\rho}(x',\cdot,\omega)|_{H^2(Y^\rho)}\le c|x-x'||Y^\rho|^{1/2}
\] 
for a constant $c$ independent of $\rho$ and $\omega$. We thus have $|u_1^\rho(\cdot,\cdot,\omega)|_{\hat{\cal H}}\le c|Y^\rho|^{1/2}$ for a constant $c$ independent of $\omega$ and $\rho$. 
%
%
\end{appendices}
\vskip 30pt
{\bf Acknowledgment} VHH and WCT gratefully acknowledge the financial support of the Singapore Ministry of Education Academic Research Fund Tier 2 grant MOE2017-T2-2-144. CHP is supported by a Graduate Scholarship of Nanyang Technological University. 

\bibliographystyle{plain}
\bibliography{paper}

\end{document}